# VASSILIEV INVARIANTS I : BRAID GROUPS AND RATIONAL HOMOTOPY THEORY

by Louis FUNAR

## CHAPTER 0

### Introduction

This paper is the first one in a series of papers on Vassiliev invariants and we are dealing here only with braid groups. The present article is a fairly detailed account starting with Chen's theory of iterated integrals and Kontsevich's approach for the universal Vassiliev invariant. We get the geometric construction of Malcev's completion (over $\mathbb{Q}$) of a discrete group in order to apply it for the case of pure braid groups $P_n$. Our first main ingredient in reconstructing Vassiliev invariants is the canonical arrow

$$P_n \longrightarrow U(P_n \otimes \mathbb{Q})$$

which we further identify with the universal Vassiliev invariant for pure braids. This is certainly transparent, even if never explicitly stated, in previous work of Stanford [St], Bar-Natan [BN2], Cartier [Car] and Kohno [Kohno1].

The extension of this morphism to the whole braid group $B_n$ cannot be a homomorphism. The reason is that $B_n$ has not a cohomology group large enough (to inject into its completion). This is equivalent to saying that only multiplicative Vassiliev invariants do not suffice to classify braids, as was the case with the pure braids. We notice however that we can build up a representation of $B_n$ related to Vassiliev invariants using the formulas of Drinfeld (see [Drin]) as is done in [Piu]. Anyway we may extend the previous morphism to a map

$$B_n \longrightarrow V(n) = U(P_n \otimes \mathbb{Q}) \ltimes S_n$$

whose failure to be a morphism might be explicitly computed in terms of Drinfeld associator.

We shall only discuss some points about regularizing singular integrals following Le and Murakami in order to explain the multiplication law. A detailed construction will be given in the second paper in this series.

In this setting the chord diagram algebras will be a sort of Malcev completions for the semi-group of knots, revealing the rational homotopic nature of Vassiliev invariants. The same conclusions were obtained by Kassel and Turaev in [KT].

The only novelty in this paper, mainly expository, is just the emphasing of this relationship which will be exploited further. We already notice that Malcev's completion has an universality property: any multiplicative universal Vassiliev (for pure braids) taking values in a graded algebra $A_n$ factors through $P_n \otimes \mathbb{Q}$. This is the case for the graded algebra $AP_n$ from [BN2] of chinese character diagrams. This means that up to an automorphism of $P_n \otimes \mathbb{Q}$ any universal invariant of pure braids has in its expansion only horizontal and Lie polynomial chord diagrams, supporting the conjecture of Bar-Natan [BN2].

*Acknowledgements.* — I would like to thanks Valentin Poénaru, Stefan Papadima, Gregor Masbaum for the stimulating discussions we had on Vassiliev invariants. The participants of the topology seminar at Grenoble were sufficiently patient to listen me, their critics and suggestions being an impulse to write this introduction to Vassiliev invariants from a rational homotopic perspective. I thanks them all. Nethertheless I'am indebted to Arlette Guttin-Lombard for typing this text.

# CHAPTER 1

## Review of Chen's theory

**1. Setup.**

(1.0)   Let $X$ be a connected $C^\infty$-manifold having $H_*(X)$ and $\pi_1 X$ finite generated. Chen ([Chen1]) constructed a series of simply-connected nilpotent Lie groups

$$\cdots \mathcal{G}(r) \longrightarrow \mathcal{G}(r-1) \longrightarrow \cdots \longrightarrow \mathcal{G}(1)$$

and a sequence of locally flat connections on $X$ which lead to holonomy homomorphisms

$$\pi_1(X) \longrightarrow \mathcal{G}(r), \quad r \geq 1.$$

If $\pi_1 X$ is torsion free nilpotent then $\mathcal{G}(r)$ stabilizes for large $r$ and the corresponding holonomy homomorphism sends $\pi_1 X$ isomorphically into $\mathcal{G}(r)$ as an uniform discrete subgroup so that $\mathcal{G}(r)$ may be identified with Malcev's completion of $\pi_1 X$.

(1.1)   If the De Rham complex $\Lambda^*(X)$ of $X$ is equipped with a direct sum decomposition of the type

$$\Lambda^p(X) = H^p \oplus dA^{p-1} \oplus A^p$$

such that $H^p$ consists in closed $p$-forms and $A^p$ contains no non-zero closed $p$-form then the locally flat connections mentioned above are uniquely determined. For a compact Riemann manifold there is a canonical decomposition of this type, namely the Hodge decomposition.

(1.2)   Let $V$ be the graded vector space with

$$V_q = H_{p+1}(X; \mathbf{k}), \quad \mathbf{k} \text{ being a fixed field},$$

and $\overline{T}(V)$ be the completion of the tensor algebra on $V$. Then every direct sum decomposition as above gives rise to a canonical differential mapping

$$\partial : \overline{T}(V) \longrightarrow \overline{T}(V)$$

having the degree -1, and a canonical $\overline{T}(V)$-valued formal power series connection $\omega$ which is a twisting cochain i.e.

$$\partial \omega + K(\omega) = 0$$

where $K$ denotes the curvature of the connection. The holonomy homomorphism is a chain map from the smooth chain complex $C_*(\Omega X)$, of the loop space $\Omega X$, to $\overline{T}(V)$, which induces an isomorphism

$$H_*(\Omega X; \mathbf{k}) \simeq H_*(\overline{T}(V))$$

in the case where $X$ is 1-connected. In the non simply-connected case there is an induced morphism

$$\mathbf{k}\pi_1 X = H_0(\Omega X; \mathbf{k}) \longrightarrow H_0(\overline{T}(V)) \stackrel{\text{not}}{=} \mathcal{U}.$$

$\mathcal{J}$ denotes the augmentation ideal of the group algebra $\mathbf{k}\pi_1 X$ and $\mathcal{U}_s$ denotes the quotient of $\mathcal{U}$ by the $s$-th power of its augmentation ideal.

THEOREM 1. — *The following sequences*
$$0 \longrightarrow \mathcal{J}^{s+1} \longrightarrow \mathbf{k}\pi_1 X \longrightarrow \mathcal{U}_s \longrightarrow 0$$
$$0 \longrightarrow \bigcap_s \mathcal{J}^s \longrightarrow \mathbf{k}\pi_1 X \longrightarrow \mathcal{U} \longrightarrow 0$$
*are exact sequences for any $s \geq 1$.*

Remark 2. — Let $G_{(r)} = \{g; g - 1 \in \mathcal{J}_G^r\}$ where $\mathcal{J}_G$ is the augmentation ideal of $\mathbf{k}G$ and let $G_*$ denotes the lower central series of the group $G$ (defined by $G_0 = G$, $G_{r+1} = [G_r, G]$). Then a general result states that
$$G_r \subseteq G_{(r)} \quad \text{for any} \quad r$$
hence $G/G_{(r)}$ is torsion free nilpotent. It is true that
$$\bigcap_r \mathcal{J}_G^r = 0 \quad \text{if and only if} \quad \bigcap_r G_r = 0$$
or, equivalently $G$ is residually torsion free nilpotent. As a consequence.

COROLLARY 3. — *If $\pi_1 X$ is residually torsion free nilpotent then the map $\mathbf{k}\pi_1 X \to \mathcal{U}$ is injective.*

## 2. Some definitions and technicalities.

(2.0) Let $V_*$ be a graded vector space, $X_1, X_2, \ldots$ be a basis for $V_*$ so that $X_1, \ldots, X_m$ form a basis for $V_0$, $X_{m+1}, \ldots, X_{m+\ell}$ form a basis for $V_1$, etc. Let $\overline{T}(V)$ be the completion of the graded tensor algebra on $V_*$; regard $X_1, X_2, \ldots$ as non-commutative variables and write $X_{i_1} X_{i_2} \cdots X_{i_p}$ for $X_{i_1} \otimes X_{i_2} \otimes \cdots \otimes X_{i_p}$. Then every element of $\overline{T}(V)$ is a formal power series $a = a_0 + \sum_i a_i X_i + \sum_{i,j} a_{ij} X_i X_j + \cdots$. The augmentation map is
$$\overline{T}(V) \longrightarrow \mathbf{k}, \quad a \longmapsto a_0$$
so
$$\mathcal{J}^r = \{a; a_{i_1 \cdots i_s} = 0 \text{ if } s < r\}.$$
We topologize $\overline{T}(V)$ using the system of neighborhoods $\{\mathcal{J}^r; r = 1, 2, \ldots\}$ of 0 so that $\overline{T}(V)$ is Hausdorff. A derivation $\partial$ of $\overline{T}(V)$ is a linear endomorphism of degree -1 satisfying the usual Leibniz rule
$$\partial(uv) = (\partial u)v + (-1)^{\deg u} u \partial v$$
and also
$$\partial \text{ is continuous and } \partial \overline{T}(V) \subset \mathcal{J}.$$

(2.1) We come back to the case when $V_*$ is the graded homology vector space of a $C^\infty$-manifold $X$. Consider the endomorphism $J : \Lambda^* X \to \Lambda^* X$ of the De Rham complex, given by
$$Jw = (-1)^{\deg w} w.$$
Let us denote by $\overline{T}_{\Lambda(X)}(V))$ the algebra of $\overline{T}(V)$-valued forms on $X$.

A formal connection on $X$ is an element $\omega \in \overline{T}_{\Lambda(X)}(V)$
$$\omega = \sum w_i W_i + \sum w_{ij} X_i X_j + \cdots$$

such that

$$w_{i_1 \cdots i_r} \text{ is a form on } X \text{ of degree } 1 + \deg X_{i_1} + \cdots + \deg X_{i_r}.$$

The curvature of the connection $\omega$ is defined as

$$K(\omega) = d\omega - J\omega \wedge \omega \in \overline{T}_{\Lambda(X)}(V)$$

where

$$d\omega = \sum dw_i X_i + \sum dw_{ij} X_i X_j + \cdots$$
$$J\omega = \sum Jw_i X_i + \sum Jw_{ij} X_i X_j + \cdots$$

Suppose a decomposition as in (1.1) is fixed:

$$\Lambda^*(X) = H^* \oplus dA^{*-1} \oplus A^*.$$

Choose a basis $X_1 X_2, \ldots$ of $V_* = H_{*+1}(X; \mathbf{k})$ and the forms $w_i$ in $\Lambda^*(X)$ so that their cohomology classes $[w_i]$ in $H^*(X)$ are representing a dual basis of $X_i$ in $H_*(X)$. We can furthermore choose $w_i$ in $H^* \subset \Lambda^*(X)$. Then the element

$$\beta = \sum w_i X_i \in \overline{T}_{\Lambda(X)}(V)$$

is independent on the choices of bases we have done, being uniquely determined by the decomposition.

THEOREM 4. — *There exists uniquely a formal connection*

$$\omega = \sum w_i X_i + \sum w_{ij} X_i X_j + \cdots$$

*and a derivation $\partial$ of $\overline{T}(V)$ such that*

(1)    *the initial term $\sum w_i X_i$ is $\beta$.*

(2)    $w_{ij}, w_{ijk}, \ldots$ *belong to $A^*$.*

(3)    $\partial \omega + K(\omega) = 0$ *(the flatness of $\omega$).*

Remark 5. — The cup-product on $H^*(X)$ determines the first stage of $\partial$ in the following way: assume that $[Jw_i \wedge w_j] = \sum c_{ij}^k [w_k]$. Then

$$\partial X_\lambda = \sum c_{ij}^\lambda X_i X_j + \cdots$$

(2.2)    Observe that $\overline{T}(V_0) = \overline{T}(V)_0$ is an (ungraded) algebra and

$$\mathcal{J}_0 = \mathcal{J} \cap \overline{T}(V)_0.$$

Let $\mathcal{N}$ be the closure of the ideal generated by $\partial V_1 \subset \overline{T}(V)_0$. Then $\mathcal{N} = \partial(\overline{T}(V)_1)$ since $\partial V_0 = 0$ and we have

$$\mathcal{U} = H_0(\overline{T}(V), \partial) = \overline{T}(V_0)/\mathcal{N}.$$

The augmentation ideal of $\mathcal{U}$, $\mathcal{J}_{\mathcal{U}} = \mathcal{J}_0/\mathcal{N}$. Then the algebras

$$\mathcal{U}_s = \mathcal{U}/\mathcal{J}_{\mathcal{U}}^{s+1}$$

are finite dimensional. Set $\nu_r : \overline{T}(V_0) \longrightarrow \mathcal{U}_r$ for the natural projections, $\mathcal{N}_{\Lambda(X)} = \Lambda^*(X) \otimes \mathcal{N} \subset \overline{T}_{\Lambda(X)}(V)$, $\nu_r : \overline{T}_{\Lambda(X)}(V) \longrightarrow \Lambda^*(X) \otimes \mathcal{U}_r$ for the natural extension of $\nu_r$. Let $L(V)$ be the graded free Lie algebra generated by $X_1, X_2, \ldots$ and $\overline{L}(V)$ its topological closure in $\overline{T}(V)$, and $\overline{L}_{\Lambda(X)}(V)$ the space of $\overline{L}(V)$-valued forms on $X$. Lets define

$$\mathcal{J}_{\mathcal{U}_r} = \nu_r(\overline{L}(V_0)) \subset \mathcal{U}_r.$$

Since $\mathcal{J}_{\mathcal{U}_r}$ is nilpotent it follows that $g_r \subset \mathcal{J}_{\mathcal{U}_r}$ is a nilpotent Lie algebra hence $\mathcal{G}(r) = \exp g_r$ is a simply-connected Lie group. This is the tower of Lie groups from (1.0).

## CHAPTER 2

## Review of Vassiliev invariants

**1. Combinatorics.**

(1.0) Any **k**-valued invariant $V$ of oriented knots in $S^3$ (or, more generally in a 3-manifold $M^3$) can be extended canonically to be an invariant of immersed circles in $S^3$ which have only ordinary double points using the following resolution of a local singularity

$$V\left(\times\!\!\times\right) = V\left(\times\!\!\times\right) - V\left(\times\!\!\times\right).$$

As usually such a skein relation means that $\times\!\!\times$, $\times\!\!\times$ or $\times\!\!\times$ are parts of bigger graphs which are identical outside a small sphere, inside of which they look as in the figures.

(1.1) Let $m \in \mathbb{Z}_+$. An invariant $V$ of oriented knots is called an invariant of type $m$ (or a Vassiliev invariant of degree $m$) if $V$ vanishes on singular knots that have more than $m$ double points:

$$V\left(\underbrace{\times\!\!\times\,\times\!\!\times\,\ldots\,\times\!\!\times}_{>m}\right) = 0.$$

The **k**-space of Vassiliev invariants $\mathcal{V}$ is the space of invariants of finite degree and has a natural filtration by the degree $\mathcal{V}_*$.

(1.2) A chord diagram is an oriented circle with finitely many chords marked on it regarded up to orientation preserving diffeomorphisms of the circle. Denote by $D$ the collection of all chord diagrams graded by the number of chords.

Now a **k**-weight system of degree $m$ is a function

$$W : D_m \longrightarrow \mathbf{k}$$

which fulfills :

(1) If $d \in D_m$ has an isolated chord (which does not intersect the other chords of $d$) then

$$W(d) = 0.$$

(2) Whenever four diagrams $d_1$, $d_2$, $d_3$, $d_4$ differ only as shown in the figure below, their weights satisfy the $4T$-relation

$$W(d_1) - W(d_2) = W(d_3) - W(d_4)$$

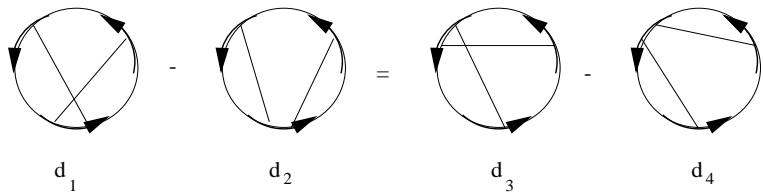

Let $\mathcal{W}_*$ denotes the graded space of weight systems, and $gr^*\mathcal{V} = (\mathcal{V}_*/\mathcal{V}_{*+1})$ be the graded **k**-space of Vassiliev invariants obtained from its natural filtration.

enables us to consider $\mathcal{J}_k(\mathcal{N})$ the ideal generated by the images of singular knots having fewer than $k$ self-crossings, under repeated use of the desingularisation. Then $\mathcal{J}_*(\mathcal{N})$ is an ascending filtration whose graded space $gr^*\mathcal{N} = (\mathcal{J}_*(\mathcal{N}))/\mathcal{J}_{*+1}(\mathcal{N}))_*$ is isomorphic with $gr^*\mathcal{V}$.

(1.4) The main theorem in Vassiliev theory is ([BN1], [Konts], [Vass], [Bir-Lin], [V]). ables us to consider $\mathcal{J}_k(\mathcal{N})$ the ideal generated by the images of singular knots having fewer than $k$ self-crossings, under repeated use of the desingularisation. Then $\mathcal{J}_*(\mathcal{N})$ is an ascending filtration whose graded space $gr^*\mathcal{N} = (\mathcal{J}_*(\mathcal{N}))/\mathcal{J}_{*+1}(\mathcal{N}))_*$ is isomorphic with $gr^*\mathcal{V}$.

(1.4) The main theorem in Vassiliev theory is ([BN1], [Konts], [Vass], [Bir-Lin]).

THEOREM 1. — *We have an isomorphism of graded $\mathbf{k}$-spaces*

$$\mathcal{W}^* \simeq gr^*\mathcal{V}$$

*for $\mathbf{k} = \mathbb{R}$.*

Actually this result was improved at $\mathbf{k} = \mathbb{Q}$ by Le-Murakami ([Le-Mu]). The main feature of Vassiliev invariants over $\mathbb{R}$ is that via Theorem 1 they are algorithmically computable using the so-called actuality tables (see [Bir], [Big-Lin]).

## 2. The algebra of diagrams.

(2.0) Set $\mathcal{D} = \mathbf{k}D$ for the $\mathbf{k}$-space spanned by chord diagrams and

$$\mathcal{A} = \mathcal{D}/\mathbf{k}(4T\text{-relations})$$

the $\mathbf{k}$-algebra of diagrams, and set for the algebra of reduced diagrams $\mathcal{A}^r$,

$$\mathcal{A}^r = \mathcal{A}/(d \text{ having isolated chords}).$$

It is clear that the weight system are actually functionals on $\mathcal{A}^r$.

(2.1) The multiplication of two diagrams in $\mathcal{D}$ is obtained by connected sum of diagrams in two points not lying on any chord. The ambiguity is cancelled when passing to the quotient $\mathcal{A}$ and gives rise to a multiplication $\mathcal{A} \times \mathcal{A} \to \mathcal{A}$.

We have a co-multiplication $\Delta : \mathcal{A} \to \mathcal{A} \otimes \mathcal{A}$ by

$$\Delta(d) = \sum d' \otimes d''$$

where $d'$ is obtained from $d$ by deleting some chords, $d''$ by deleting the chords of $d'$ and the sum being taken over all possibilities.

Remark that a natural multiplication and co-multiplication may be defined in a similar vein for $\mathcal{N}$.

(2.2) The interest in having much structure on $\mathcal{A}$ is to derive a simpler algebraic description of it. In fact

THEOREM 2. — $(\mathcal{A}, \cdot, \Delta)$ *is a commutative and co-commutative Hopf algebra over* $\mathbf{k}$.

Therefore by the structure of Hopf algebras we know that $\mathcal{A}$ is the symmetric algebra generated by the primitive elements of $\mathcal{A}$:

$$\mathcal{A} = S(P(\mathcal{A})), \quad P(\mathcal{A}) = \{a \in \mathcal{A};\ \Delta(a) = a \otimes 1 + 1 \otimes a\}.$$

If $\mathcal{A}^*$ is the dual Hopf algebra and $P'(\mathcal{A}^*)$ is the set of primitive elements of degree greater than 1 then we can identify the $\mathbf{k}$-space of weight systems as

$$\mathcal{W} \approx S(P'(\mathcal{A}^*)).$$

Notice that the above isomorphism is a graded isomorphism, and for $\mathbf{k} = \mathbb{R}$ using Theorem 1, the weight systems in $P'(\mathcal{A}^*)$ correspond to Vassiliev invariants which are additive under the operation of taking connected sum of knots.

(2.3) A computer search gave for the number $d_n$ of primitive elements in degree $n$ the values (see [BN1]):

$$d_1 = 1,\ d_2 = 1,\ d_3 = 1,\ d_4 = 2,\ d_5 = 3,\ d_6 = 5,\ d_7 = 8,\ d_8 = 12.$$

We set $P_{\mathcal{A},\mathbf{k}}$ for the Hilbert series of $\mathcal{A}$, which is

$$P_{\mathcal{A},\mathbf{k}}(t) = \prod_{i=1}^{\infty}(1 + t^{d_i})$$

because $\mathcal{A}$ is a polynomial algebra.

### 3. Kontsevich's universal invariant.

(3.0) We outline below the construction of the isomorphism of Theorem 1. The easy part is to start with a Vassiliev invariant of degree $m$, say $V$ and to derive a weight system.

Let $d \in D_m$ be a chord diagram. an embedding of in $\mathbb{R}^3$ is an immersion $i_d : S' \to \mathbb{R}^3$ whose singularities are ordinary double points and satisfies:

$i_d(a) = i_d(b)$ iff $a = b$ or else $a$ and $b$ are the endpoints of a chord in $d$.

There exists an unique regular homotopy class of such immersions for a fixed chord diagram hence any two embeddings $i_d$ and $\tilde{i}_d$ are connected by a sequence of flips in which an over-crossing 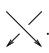 changes into an under-crossing.

(3.1) *Example*:

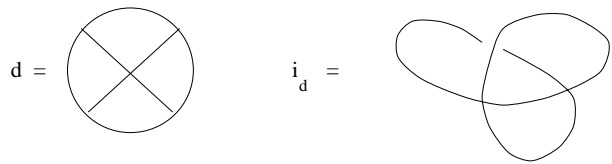

(3.2) Let us define
$$w(d) = V(i_d), \text{ for } d \in D_m, V \in \mathcal{V}_m.$$

A flip does not change the value of $V(i_d)$ since

$V(i_d) - V(\tilde{i}_d) = V$ ( a singular knot with $m+1$ double points) $= 0$,

so $V(i_d)$ is well-defined. It remains to see that $w$ is actually a weight system.

Firstly we have
$$V\left(\begin{array}{c}\text{\fbox{A}}\bowtie\text{\fbox{B}}\end{array}\right) = 0$$

because 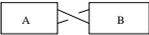 and 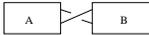 are isotopic, which implies that $w(d) = 0$ if $d$ contains an isolated chord.

Further consider $K_0$ and $K_1$ be two knots with $m-1$ double points which are identical outside a small sphere, inside which they look as in figure below:

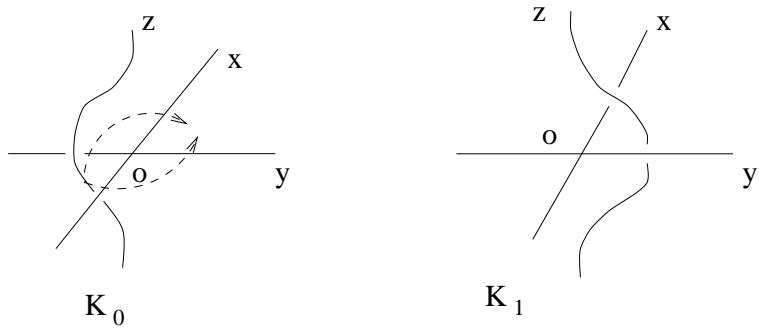

We can pass from $K_0$ to $K_1$ in two ways: by moving the strand $z$ to cross $x$ and then $y$, or to cross $y$ firstly and further $x$. But each time $z$ cross $x$ or $y$ we can compute the change in $V(K_0)$ using the values of $V$ on knots with $n$ double points. The two ways to get $V(K_1) - V(K_0)$ must give the same answer hence we derive a four term relation on $w$ which is exactly $(4T)$.

(3.3) The inverse homomorphism $\mathcal{W} \to \mathcal{V}$ is provided by using the Knizhnik-Zamolodchikov equation and is due to Kontsevich. It gives a sort of universal link invariant taking values in $\mathcal{W}$.

(3.4) Recall that Chen [Chen2] gives an effective method to compute the holonomy of a flat connection $\Omega$ on a $C^\infty$-manifold $X$, taking values in a topological algebra $A$ over $\mathbf{k}$, with unit 1. The parallel transport along the smooth curve $\gamma : I \to X$ is the map
$$h_\Omega(\gamma) : I \longrightarrow A, \quad I = [0,1]$$
which satisfies
$$h_\Omega(0) = 1, \quad \frac{\partial}{\partial t}h_\Omega(t) = \Omega(\dot{\gamma}(t))h_\Omega(t), \quad t \in I$$
if such a function exists and it is unique.

If $\Omega$ is flat then this holonomy map $h_\Omega(\gamma)$ is invariant under homotopies of $\gamma$ which preserve its endpoints, and it can be calculated by means of iterated path integrals as:
$$h_\Omega(\gamma) = 1 + \sum_{m=1}^{\infty} \int_{0 \leq t_1 < t_2 < \cdots < t_m \leq 1} \gamma^*\Omega(t_1) \wedge \cdots \wedge \gamma^*\Omega(t_m)$$
where $\gamma^*\Omega(t_1) \wedge \cdots \wedge \gamma^*\Omega(t_m)$ is a top form on the simplex $\Delta_m$.

(3.5) Let $X_n$ be the configuration space of $n$ distinct points in $\mathbb{C}$,

$$X_n = \{(z_1, \ldots, z_n) \in \mathbb{C}^n, \ z_i \neq z_j, \ \forall i \neq j\} \subset \mathbb{C}^n$$

and $\omega_{ij} \in \Lambda^1(X_n)$ defined by

$$\omega_{ij} = \frac{dz_i - dz_j}{z_i - z_j}.$$

Let $D_n^{KZ}$ be the collection of all diagrams made by $n$ ordered downward pointing arrows and arcs connecting them (with eventually 3-valent cyclically orientations around vertices):

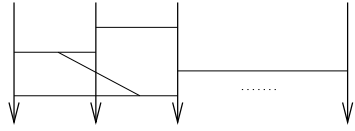

and set

$$\mathcal{A}_n^{KZ} = \mathbf{k}(D_n^{KZ})/(STU\text{-relations})$$

where the $STU$ relations are

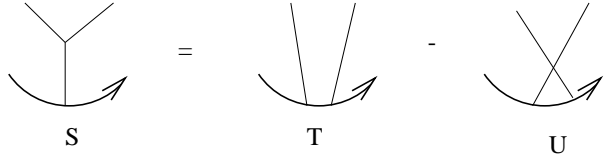

S        T        U

This definition is reminiscent to the identification of $\mathcal{A}$ with the algebra of 3-valent diagrams $\mathcal{A}^t$ where the $4T$-relation is replaced by the $STU$-relation, due to Kontsevich and Bar-Natan. We refer to [BN1] for more details. Remark only that $\mathcal{A}_n^{KZ}$ has a multiplication by putting diagrams one above the other.

Set

$$\Omega_{ij} = \Big\downarrow \ \Big\downarrow \ \cdots \ \underset{i}{\Big\downarrow}\!\!\overline{\phantom{\cdots}}\!\!\underset{j}{\Big\downarrow} \ \cdots \ \Big\downarrow \ \in \mathcal{A}_n^{KZ}$$

and form the formal Knizhnik-Zamolodchikov connection on $X_n$

$$\Omega_n = \sum_{1 \leq i < j \leq n} \Omega_{ij}\omega_{ij} \in \Lambda^1(X_n) \otimes \mathcal{A}_n^{KZ}.$$

Then the $STU$-relation gives the flatness of $\Omega_n$.

(3.6) The connection $\Omega_n$ has a simple generalization to the case when the underlying algebra is $\mathcal{A}_{n,n}^{KZ}$ generated by diagrams having $2n$ arrows whose $n$ arrows point upward and whose remaining $n$ arrows point downward. Then one defines

$$\Omega_{n,n} = \sum_{i<j} s_i s_j \Omega_{ij} \omega_{ij} \in \mathcal{A}_{n,n}^{KZ} \otimes \Lambda^1(X_n)$$

where $s_i = \begin{cases} 1 & \text{if the } i^{\text{th}} \text{ arrow points downward} \\ -1 & \text{otherwise} \end{cases}$.

In general we can specify the signature of arrows as $\varepsilon : \{1, \ldots, 2n\} \to \{\pm 1\}$ and identifying $\mathcal{A}_{n,n}^{KZ}$ with a specific $\mathcal{A}_{n,n}^{KZ}(\varepsilon)$.

(3.7) Choose a decomposition $\mathbb{R}^3 = \mathbb{C} \times \mathbb{R}$ and let $K : S^1 \to \mathbb{R}^3$ be a parametrized knot whose projection on $\mathbb{R}$ is a Morse function. Consider the series

$$Z(K) = \sum_{m=0}^{\infty} \frac{1}{(2\pi i)^m} \int_{t_{\min} \leq t_1 < \cdots < t_m \leq t_{\max}} \sum_{P \text{ pairing } \{(z_i, z'_i)\}} (-1)^{\#P \downarrow} D_P \bigwedge_{i=1}^{m} \frac{dz_i - dz'_i}{z_i - z'_i} \in \mathcal{A}^r_{\mathbb{C}}$$

where

(1) the projection of $K$ on $\mathbb{R}$ is $[t_{\min}, t_{\max}]$;

(2) a pairing $P$ is a choice of unordered pairs $(z_i, z'_i)$, $1 \leq i \leq m$ for which $(z_i, t_i)$ and $(z'_i, t_i)$ are distinct points of $K$;

(3) $\#P \downarrow$ is the number of points in the pairing $P$ where the orientation of $K$ points downward with respect to the projection on $\mathbb{R}$;

(4) $\mathcal{A}^r_{\mathbb{C}}$ is $\mathcal{A}^r$ for $\mathbf{k} = \mathbb{C}$;

(5) $D_P$ is the diagram associated to the $m$ pairs of points in $S^1$.

We defer to [BN1] for a complete proof that $Z_n(K)$ is well-defined and invariant ot homotopies which preserve the number of critical points. Let the symbol $\infty$ states for the embedding and notice that $Z(\infty) = 0+$ (higher order terms) hence $Z(\infty)$ is invertible in $\mathcal{A}^r_{\mathbb{C}}$. For $K$ an embedded Morse knot with $c$ critical points we set

$$\widetilde{Z}(K) = Z(\infty)^{1-\frac{c}{2}} Z(K) \in \mathcal{A}^r_{\mathbb{C}}$$

which is an isotopy invariant, called the universal Kontsevich invariant. It is simply to check that $\widetilde{Z}$ provides the inverse morphism $\mathcal{W}^* \to gr^* \mathcal{V}$. Observe that $t \to -t$, $z \to \bar{z}$ maps a knot into an equivalent one while $\Omega_{n,n} \to -\overline{\Omega}_{n,n}$. This proves that $\widetilde{Z}(K) \in \mathcal{A}^r_{\mathbb{R}}$.

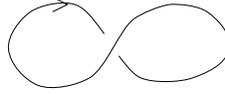

# CHAPTER 3

## The Malcev completion of the group of pure braids

### 1. Configuration spaces.

(1.0) We come back to Chen's theory from the first chapter. We want to carry out this general theory in the specific case of configuration spaces. It is known that

$$\pi_1 X_n = P_n$$

is the group of pure braids in $n$ strings. Its cohomology ring was computed by Arnold ([Arnold]) and it is

$$H^*(P_n; \mathbb{Z}) = \langle e_{ij}, \ 1 \leq i < j \leq n, \ \deg e_{ij} = 1; \ e_{ij} e_{jk} = e_{ik} e_{jk} + e_{ij} e_{ik}, \ 1 \leq i < j < k \leq n \rangle.$$

Under the natural map $H^*(P_n; \mathbb{Z}) \hookrightarrow H^*(X_n, \mathbb{C}) \approx H^*_{DR}(X_n)$ the generators $e_{ij}$ correspond to the 1-forms $\omega_{ij} \in \Lambda^1(X_n)$.

From this description, or directly by using the fact that $P_n$ is an amalgamation of free groups $\mathbb{F}_{n-1} \ltimes \mathbb{F}_{n-2} \ltimes \cdots \ltimes \mathbb{F}_1$ we derive that the Hilbert polynomial is

$$P_{X_n, \mathbf{k}}(t) = (1+t)(1+2t) \cdots (1+(n-1)t)$$

over any field $\mathbf{k}$.

(1.1) Let $A_{n,\ell} = \{\{(i_1, j_1), (i_2, j_2), \ldots, (i_\ell, j_\ell)\}\}$ lexicographically ordered, to avoid permutations, where $i_s, j_s \in \{1, \ldots, n\}$ satisfy

(1) $\quad i_s < j_s$;

(2) $\quad \{i_1, \ldots, i_\ell\} \cap \{j_1, \ldots, j_\ell\} = \emptyset$.

For a multi-index $I \in A_{n,\ell}$ we set $e_I = e_{i_1 j_1} e_{i_2 j_2} \cdots e_{i_\ell j_\ell}$. It follows easily that $\{e_I; \ I \in A_{n,\ell}\}$ is a **k**-basis for $H^\ell(X_n; \mathbb{Z})$. Our first task is to define a flat formal connection $(\omega, \partial)$ on $X_n$.

We set then
$$\omega = \sum_I w_I Z_I$$

where $w_I$ is the form $\bigwedge_{s=1}^\ell \omega_{i_s j_s}$ dual to $e_I$, and $Z_I$ are formal non-commutative indeterminates with $\deg Z_I = |I| - 1$.

PROPOSITION 1. — *There exists an unique formal homological connection $(\omega, \partial)$.*

*Proof.* — The quadratic formal connection is given by the cohomology cup-product. We have (see Remark 5, chap. 1)
$$\partial Z_I = \sum c_I^{JK} Z_J Z_K$$

where $[Jw_J \wedge w_K] = \sum c_I^{JK}[w_I]$ in $H^*(X_n)$. Also we can write
$$J\omega = \sum (-1)^{|I|-1} w_I Z_I$$
$$J\omega \wedge \omega = \sum (-1)^{|I|-1} w_I \wedge w_J Z_I Z_J$$
$$d\omega = \sum dw_I Z_I = 0 \quad \text{since the forms } w_I \text{ are closed}$$
$$\partial \omega = \sum_I w_I \left( \sum_{J,K} c_I^{JK} Z_J Z_K \right) = \sum_{J,K} \left( \sum_I c_I^{JK} w_I \right) Z_J Z_K$$
$$= \sum_{J,K} (-1)^{|J|-1} w_J \wedge w_K Z_J Z_K.$$

The last equality follows from the very pleasant fact that $c_I^{JK}$ are determined for $X_n$ directly at the form level *i.e.*
$$Jw_J \wedge w_k = \sum c_{JK}^I w_I, \quad I, J, K \in \bigcup_\ell A_{n,\ell}.$$

We derive
$$d\pi + \partial \omega - J\omega \wedge \omega = 0$$

hence the flatness of $(\pi, \partial)$. ∎

Remark this proposition is equivalent to the fact that the spaces $X_n$ are formal (see [DGMS]).

(1.2) Now, with the notations of Chapter 1, the degree 0 component of $\overline{T}(V)$ is therefore a quotient of $\mathbf{k}[[Z_{(ij)}, \ 1 \leq i < j \leq n]]$ (the double brackets states for the series in non-commutative variables), and it remains to compute effectively $\partial$ on $Z_{(ij)(k\ell)}$ for obtaining $H_0(\overline{T}(V), \partial)$. These computations are giving in fact the Lie algebra of $P_n$ (see [Kohno2], [FR]) but we carry them out for the sake of completeness. It is immediate that
$$c_{\{(i,j),(k,\ell)\}}^{\{(u,v)\}\{(r,s)\}} = 0, \quad \text{if } \{i,j\} \cap \{k,\ell\} = \emptyset$$
$$c_{\{(i,j),(k,\ell)\}}^{\{(u,v)\}\{(r,s)\}} = \delta_{ik}^{uv} \delta_{jk}^{rs} + \delta_{ij}^{uv} \delta_{ik}^{rs}$$

and $c_I^- = -c_I^+$. We obtain then
$$\partial Z_{(i,j)(k,\ell)} = Z_{ij}Z_{k\ell} - Z_{k\ell}Z_{ij} \quad \text{if} \quad \{i,j\} \cap \{k,\ell\} = \emptyset$$
$$\partial Z_{(i,j)(i,k)} = Z_{ij}Z_{jk} + Z_{ij}Z_{ik} - Z_{jk}Z_{ij} - Z_{ik}Z_{ij}$$
$$\partial Z_{(i,k)(j,k)} = Z_{ij}Z_{jk} + Z_{ik}Z_{jk} - Z_{jk}Z_{ij} - Z_{jk}Z_{ik}.$$

Then the universal algebra $\mathcal{U}(n) = H_0(\overline{T}(V), \partial)$ can be presented as
$$\mathcal{U}(n) = \mathbf{k}[[Z_{ij},\ 1\leq i<j\leq n]]/[Z_{ij}, Z_{k\ell}] = 0 \quad \text{if} \quad \{i,j\} \cap \{k,\ell\} = \emptyset$$
$$[Z_{ij} + Z_{jk}, Z_{ik}] = 0 \quad \text{if} \quad 1\leq i<j<k\leq n$$
$$[Z_{ik} + Z_{ij}, Z_{jk}] = 0.$$

Its augmentation ideal $\mathcal{J}_{\mathcal{U}(n)} = (Z_{ij};\ 1\leq i<j\leq n) \subset \mathcal{U}(n)$.

(1.3) Remark that this algebra is well-understood object. Again we can use fairly general results of Kohno, Falk & Randell, see also Berceanu ([Kohno2], [FR], [Ber]), or else to use the semi-direct product decomposition of $P_n$ into free groups (and to use that at algebra level the semi-direct product transforms into direct sum) to derive the Hilbert series of the graded algebra $\mathcal{U}(n)$ as
$$P_{\mathcal{U}(n),\mathbf{k}}(t) = \prod_{j=1}^{n-1} \frac{1}{j - jt}.$$

Remark that the graduate structure on $\mathcal{U}(n)$ is given by the degree of a polynomial in $Z_{ij}$'s. This means that we computed actually the series of $gr^*\mathcal{U}(n) = \bigoplus_r \mathcal{J}^r_{\mathcal{U}(n)}/\mathcal{J}^{r+1}_{\mathcal{U}(n)}$.

Observe that $P_{\mathcal{U}(n),\mathbf{k}}$ does not depend on the characteristic of $\mathbf{k}$. If we set $\mathcal{U}(n)_\mathbb{Z}$ for the algebra defined over $\mathbb{Z}$ we deduce (see [Ber]):

COROLLARY 2. — *The algebra $\mathcal{U}(n)_\mathbb{Z}$ is torsion-free, or equivalently, the graduation $gr^*P_n = \bigoplus (P_n)_{(r)}/(P_n)_{(r+1)}$ is torsion-free.*

We have also a holonomy homomorphism
$$Z : P_n \longrightarrow \mathcal{U}(n).$$

COROLLARY 3.

1) *The holonomy $Z$ is injective.*

2) *Let consider*
$$Z\nu_r : P_n \longrightarrow \mathcal{U}(n)_r = \mathcal{U}(n)/\mathcal{J}^{r+1}_{\mathcal{U}(n)}.$$

*Then* $\ker Z\nu_r = (P_n)_r = (P_n)_{(r)}$, *(see [FR2], [St]).*

This follows directly from Chen's work and the fact that $P_n$ is a residually torsion free nilpotent group.

We can be more precise on the image of $P_n$ under $Z$. Let $g(n)$ be the free Lie algebra on $Z_{ij}$ quotiened by the ideal defining $\mathcal{U}(n)$. We can view $g(n) \subset \mathcal{U}(n)$ and $\mathcal{U}(n)$ is identified this way with the enveloping algebra of $g(n)$. Then $\mathcal{G}(n) = \exp g(n) \subset \mathcal{U}(n)$ is the set

closure of $\{\exp(x);\ x \in g(n)\} \subset \mathcal{U}(n)$

endowed with the multiplication induced by Campbell-Hausdorff formula. Then, according to Chen $P_n$ is an uniform discrete subgroup of the infinite dimensional formal Lie group $\mathcal{G}(n)$. Remark that we may define analogously $\mathcal{G}(n)_r \subset \mathcal{U}(n)_r$ which are nilpotent simply-connected Lie groups and $\mathcal{G}(n) = \varprojlim \mathcal{G}(n)_r$ so it inherits a natural topology as a closed Lie group. When $\mathbf{k} = \mathbb{Q}$, $g(n)$ is the Malcev Lie algebra of $P_n$ (see [Kohno2]).

## 2. Vassiliev invariants for pure braids.

(2.0) We can regard pure braids as isotopy classes of $n$ strings in $\mathbb{R}^2 \times [0,1]$ having nowhere horizontal tangent vectors and fixed endpoints. Therefore every invariant for braids $V$ can be uniquely extended to singular braids (where several self-crossings which are ordinary double points are allowed) by means of the skein relation we encountered in Chapter 2, namely

$$V\left(\times\!\!\!\times\right) = V\left(\times\right) - V\left(\times\right),$$

where we supposed all strings are oriented downward for the moment. Let consider a similar theory as in Chapter 2 for pure braids, so define $\mathcal{V}^*(P_n)$ the **k**-space of Vassiliev invariants of finite type, where $\mathcal{V}^m$ is spanned by those invariants vanishing on singular braids having more than $m$ double points.

(2.1) There is a natural candidate for the analogue diagrams: we denote by $\mathcal{D}_n = \mathcal{D}(\underbrace{|||\cdots|}_{n})$ the span of diagrams consisting in $n$ arrows which are labeled $1, \ldots, n$ and points downward and several horizontal arcs whose endpoints are on the arrows. We can multiply these diagrams by putting them one above the other. Consider further the algebra

$$\mathcal{T}_n = \mathcal{T}\left(\underbrace{\downarrow\downarrow\cdots\downarrow}_{n}\right) = \mathcal{D}_n / \mathbf{k}\text{ (relations (1) and (2))}$$

where

[diagram: horizontal arcs on strands $i,j,k,l$ showing commutation relation]

so horizontal arcs having distinct sets of endpoints commute, and

[diagram: six-term relation on strands $i,j,k$]

In a similar vein we can define the **k**-space of weight systems $\mathcal{W}^*(P_n)$ as functionals on $\mathcal{D}_n$ which pass to the quotient $\mathcal{T}_n$. Notice the natural grading on $\mathcal{T}_n$ is that given by the number of horizontal arcs.

PROPOSITION 4. — *We have an isomorphism $\mathcal{W}^*(P_n) \approx gr^*\mathcal{V}^*(P_n)$.*

It suffices to reread the proof of Theorem 1, chap. 2 to see that all constructions can be carried out in this simpler setting.

(2.2) PROPOSITION 5. — *We have an isomorphism of graded **k**-algebras*

$$\mathcal{U}(n) \approx \mathcal{T}_n.$$

*Proof.* — It suffices to observe that the map

$$Z_{ij} \in \mathcal{U}(n) \longmapsto \downarrow \cdots \underset{i}{\vphantom{|}}\overline{\phantom{\cdots}}\underset{j}{\vphantom{|}} \cdots \downarrow \quad \in \mathcal{T}_n$$

It remains to identify now the universal Kontsevich invariant in this setting. A comparison of the considered flat connections gives:

COROLLARY 6.

1) The homomorphism $Z: P_n \to \mathcal{U}(n)$ is the universal Kontsevich invariant for pure braids.

2) There is a 1:1 correspondence between multiplicative Vassiliev invariants (i.e. $V(xy) = V(x)V(y)$) and morphisms
$$f \in \mathrm{Hom}(P_n, \mathbf{k})$$
satisfying the property:
$$f \text{ factors as } P_n \xrightarrow{Z\nu_r} \mathcal{J}^r_{\mathcal{U}(n)}/\mathcal{J}^{r+1}_{\mathcal{U}(n)} \longrightarrow \mathbf{k}.$$

It follows from 1.3 that

COROLLARY 7 (Stanford).

1) Even for $\mathbf{k} = \mathbb{Z}$ multiplicative Vassiliev invariants classify pure braids.

2) $V(x) = V(y)$ for all multiplicative Vassiliev invariants of degree $m$ if and only if $xy^{-1} \in (P_n)_{(r)}$, or, equivalently $x \equiv y$ in the quotient $P_n/(P_n)_{(r)}$ of the lower central series.

Observe that the invariant $Z(K)$ of Morse link behaves multiplicatively with respect to connected sums. Let us be more precise: assume that in the interval $[t - \varepsilon, t + \varepsilon]$ of values there are not critical values. We tie the knot $K$ by a plane $\mathbb{R}^2 \times t$

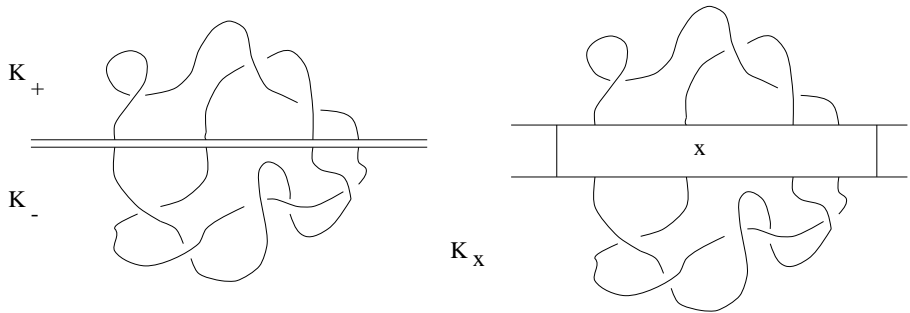

and insert a box where the strands (up and bottom ones) are connected using some pure braid $x$. We say that $K$ was modified by $x$. We have Vassiliev invariants for tangles (see [Le-Mu], [BN2]) and if we denote $K^+$ and $K^-$ the two tangles in which $K$ is splitted by $\mathbb{R}^2 \times t$ we have a multiplicativity
$$Z(K) = Z(K_+)Z(K_-)$$
$$Z(K^x) = Z(K_+)Z(x)Z(K_-)$$
where $K^x$ is the new obtained knot. Notice that the various $Z(K_{+-})$ lie in different algebras and multiplication has sense when restricted on a quotient of both. Since we can write
$$Z(K) = \sum_{m=0}^{\infty} Z_m(K)$$
where $Z_m(K)$ is the universal Vassiliev invariant of degree $m$, we see that
$$Z(x) = 1 + Z_{m+1}(x) + Z_{m+2}(x) + \cdots, \text{ if } x \in (P_n)_{(m)}.$$
Therefore we obtain

COROLLARY 8 ([St]). — Assume we modified the knot $K$ by some pure braid $x \in (P_n)_{(m)}$. Therefore all Vassiliev invariants of degree less than $m+1$ of $K$ and $K^x$ coincide.

(2.3) We have natural product operations which are a sort of exterior composition laws
$$P_n \times P_k \longrightarrow P_{n+k}$$
obtained by simply putting the strands together. On the other hand we have also
$$\mathcal{U}(n) \times \mathcal{U}(k) \longrightarrow \mathcal{U}(n+k)$$
given by:
$$(Z_{ij}, Z_{uv}) \longrightarrow Z_{ij} Z_{u+n\ v+n}.$$
With regard to the isomorphism of Proposition 5 this amounts to put then vertical arrows on the right of the $n^{\text{th}}$ arrow of the first element.

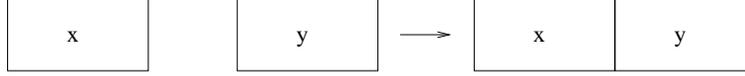

Remember we have natural injections $P_n \to \mathcal{U}(n)$, $P_k \to \mathcal{U}(k)$ which are the universal Kontsevich-Vassiliev invariants for pure braids.

PROPOSITION 9 (Product formula). — *We have a commutative diagram*
$$\begin{array}{ccc} P_n \times P_k & \longrightarrow & P_{n+k} \\ \downarrow & & \downarrow \\ \mathcal{U}(n) \otimes \mathcal{U}(k) & \longrightarrow & \mathcal{U}(n+k) \end{array}$$

*Proof.* — This follows immediately from the functoriality of Malcev completion and the easy fact that the Malcev completion of the product of two groups is the product of their Malcev completions.

Let's provide a simple geometric proof also. We seen that $Z$ comes as a monodromy representation of a flat bundle over configuration spaces. We have also $X_{n+k} \supset X_n \times X_k$. If $\Omega_{n+k}$, $\Omega_n$, $\Omega_k$ are the corresponding flat connections we have
$$\Omega_{n+k} = \Omega_n \oplus \Omega_k + \Omega^{\perp},$$
where
$$\Omega^{\perp} = \sum_{i=1}^{n} \sum_{\alpha=n+1}^{k} \Omega_{i\alpha} d\log z_i - z_{\alpha}.$$
Let consider $\gamma : [0,1] \to X_{n+k}$ which is the composition of $\gamma_1 : [0,1] \to X_n$ and $\gamma_2 : [0,1] \to X_k$. We have
$$\gamma^* \Omega_{n+k} = \gamma_1^* \Omega_n \oplus \gamma_2^* \Omega_k + \gamma^* \Omega^{\perp}.$$
Finally the $r^{\text{th}}$ iterated integral of Chen reads
$$\int_{0 \leq t_1 < \cdots < t_r \leq 1} \gamma^* \Omega_{n+k}(t_1) \wedge \cdots \wedge \gamma^* \Omega_{n+k}(t_r) = \sum_{\ell=0}^{r} \int_{0 \leq t_1 < \cdots < t_\ell \leq 1} \gamma_1^* \Omega_n(t_1) \wedge \cdots$$
$$\wedge \gamma_1^* \Omega_n(t_\ell) \int_{0 \leq t_1 < \cdots < t_{r-\ell} \leq 1} \gamma_2^* \Omega_k(t_1) \wedge \cdots \wedge \gamma_2^* \Omega_k(t_{r-\ell}) + \text{ (integrals containing } \gamma^* \Omega^{\perp}).$$
We claim that each integral from above containing $\gamma^* \Omega^{\perp}$ vanishes.

In fact the local picture of $\gamma$ contains the strands $i^{\text{th}}$ and $s^{\text{th}}$ which are not braided. We deform $\gamma$ to $\gamma_\varepsilon$ as in the picture below.

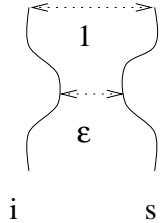

containing $s$ times $\gamma^*\Omega^\perp$ has its modulus $\sim$ constant $\frac{1}{\varepsilon^s}$. Taking $\varepsilon$ approach 0 we find the constant be 0. Notice that in the deformation $\gamma_\varepsilon$ the tangent vectors point downward. This ends the proof of Proposition 9. ∎

COROLLARY 10. — *We can compute some values of $Z : P_n \to \mathcal{U}(n)$ as:*

$$Z(b_i^2) = \exp(Z_{i\ i+1}) \in \mathcal{G}(n) \subset \mathcal{U}(n).$$

*Proof.* — It follows from the easy calculation of $Z(b_1^2)$ in $\mathcal{U}(2)$ and Proposition 9. ∎

### 3. Vassiliev invariants for braids.

(3.0) One can define the Vassiliev (or finite type) invariants for ordinary braids by requiring that their extensions to singular braids vanish on singular braids having more than $m+1$ double points, where $m$ is the degree of the invariant. Since the pure braids are distinguished by finite type invariants it follows by a straightforward argument that ordinary braids are also classified.

(3.1) However there is an important difference in the case of ordinary braids. The same approach as in the previous section fails because $B_n$ has not a cohomology ring large enough hence the holonomy morphism provided by Chen's theory is not injective. This may be rephrased by saying that multiplicative Vassiliev invariants do not classify ordinary braids.

(3.2) We shall use however the previous construction for defining a larger algebra $V(B_n)$ – which we call the Vassiliev algebra for $B_n$ – as the crossed product $\mathcal{U}(n) \ltimes S_n$, where $S_n$ is the group of permutations in $n$ letters. The extension of the homomorphism $Z : P_n \to \mathcal{U}(n)$ is a mapping $Z : B_n \to V(B_n)$ whose failure to be a group representation may be described by a sort of 2-cocycle. This 2-cocycle furnish not a representation of a group extension of $B_n$ but one of a groupoid extension.

(3.3) Remark first that $S_n$ acts on $X_n$ by permutations of the coordinates and $\pi_1(X_n/S_n) = B_n$. Unfortunately the universal flat connection $\Omega$ is not $S_n$-invariant.

On the other hand the $S_n$-action in the homology of $X_n$ induces a $S_n$-action at the tensor algebra $T(H_*(X_n))$ level. Specifically this action pass to the quotient $\mathcal{U}(n)$ and one may compute

$$S_n \times F^p \mathcal{U}(n) \longrightarrow F^p \mathcal{U}(n)$$

$$(\sigma, Z_{i_1 j_1} Z_{i_2 j_2} \cdots Z_{i_p j_p}) \longrightarrow Z_{\sigma(i_1)\sigma(j_1)} Z_{\sigma(i_2)\sigma(j_2)} \cdots Z_{\sigma(i_p)\sigma(j_p)},$$

where $F^*$ is the graduation by degree of $\mathcal{U}(n)$. This provides $\mathcal{U}(n)$ with a $S_n$-module structure.

We define now the Vassiliev algebra

$$V(B_n) = \mathcal{U}(n) \ltimes S_n.$$

Specifically as **k**-module $V(B_n)$ is $\mathcal{U}(n) \otimes \mathbf{k}[S_n]$. The product is given by

$$(Z_I, \sigma) \cdot (Z_J, \tau) = (Z_I \ ^\sigma(Z_J), \tau\sigma)$$

for multi-indices $I$ and $J$.

y. The collection $\{P(x,y), x, y \in X_n\}$ forms the fundamental groupoid of $X_n$, being endowed with a multiplication map
$$P(x,y) \times P(y,z) \longrightarrow P(x,z)$$
and a "taking the inverse" map
$$P(x,y) \longrightarrow P(y,x).$$
The parallel transport induced by the flat connection $\Omega$ furnish a series of mappings
$$Z_{xy} : P(x,y) \longrightarrow \mathcal{U}(n)$$
which is an anti-representation of the fundamental groupoid. This means that
$$Z_{xz}(uv) = Z_{xy}(u)Z_{yz}(v) \quad \text{if} \quad u \in P(x,y),\ v \in P(y,z)$$
$$Z_{ux}(u^{-1})Z_{xy}(u) = Z_{xy}(u)Z_{yx}(u^{-1}) = 1 \in \mathcal{U}(n).$$
It is clear that $Z_{xy}$ may be computed also by iterated integrals.

(3.5) We describe now the first variant to derive the universal Kontsevich-Vassiliev invariant $Z : B_n \to V(B_n)$.

Consider $z \in X_n$ be a fixed base point and $(\sigma, x) \to {}^\sigma x$ denote the $S_n$-action on $X_n$. We consider the fundamental groupoid with the set of base points $S_n z$, i.e.
$$G = \bigcup_{\sigma, \tau \in S_n} P({}^\sigma z,\ {}^\tau z).$$
Then $G$ contains several copies of $B_n$. Recall that we have an exact sequence
$$0 \longrightarrow P_n \longrightarrow B_n \overset{\sigma}{\longrightarrow} S_n \longrightarrow 0.$$
Choose some $u \in B_n$ and $\gamma$ some loop in $X_n/S_n$ representing $u$ based on $S_n z \in X_n/S_n$. Then there is an unique up to homotopy lift of $\gamma$ to a curve $\tilde{\gamma}$ in $X_n$ joining $z$ and ${}^{\sigma(u)}z$ so defining an element $\bar{u} \in P(z,\ {}^{\sigma(u)}z)$. We define furthermore
$$Z(u) = \left(Z_{z,\ \sigma(u)z}(\bar{u}), \sigma(u)\right) \in V(B_n).$$

(3.6) We have natural morphisms $B_n \times B_k \to B_{n+k}$ obtained by putting the strands together and a morphism $V(B_n) \times V(B_k) \to V(B_{n+k})$ extending the corresponding morphisms at the $\mathcal{U}(n)$-level.

PROPOSITION 11 (Product formula for $B_n$). — *We have a commutative diagram*
$$\begin{array}{ccc} B_n \times B_k & \longrightarrow & B_{n+k} \\ {\scriptstyle Z \times Z}\downarrow & & \downarrow {\scriptstyle Z} \\ V(B_n) \times V(B_k) & \longrightarrow & V(B_{n+k}). \end{array}$$

*Proof.* — The geometric proof in Proposition 9 works as well in this setting. ∎

(3.7) For a complete description of $Z$ we need to know its behavior with respect to the multiplication in $B_n$. A first step towards this is provided by:

PROPOSITION 12. — *Assume that $\sigma(u)\sigma(v) = \sigma(v)\sigma(u)$. Then we have*
$$Z(uv) = Z(u)Z(v).$$

*Proof.* — We choose $\gamma_1$ and $\gamma_2$ two loops in $X_n/S_n$ representing $u$ and $v$. Let $\tilde{\gamma}_1$ and $\tilde{\gamma}_2$ be their respective lifts in $X_n$ starting at $z$. Let $\gamma_3$ be the lift of $\gamma_2$ starting at $^{\sigma(u)}z$. Therefore $\tilde{\gamma}_1\gamma_3$ is a lift of $\gamma_1\gamma_2$ starting at $z$ and we have

$$Z(uv) = \left(Z_{z\ \sigma(vu)z}(\tilde{\gamma}_1\gamma_3), \sigma(v)\sigma(u)\right) = \left(Z_{z\ \sigma(u)z}(\tilde{\gamma}_1), \sigma(u)\right)\left(Z_{\sigma(u)z\ \sigma(vu)z}(\gamma_3), \sigma(v)\right).$$

Remark that we have an induced $S_n$-action on the fundamental groupoid $\sigma: P(x,y) \to P(^\sigma x, ^\sigma y)$. From the homotopy uniqueness of the lift we derive that $\sigma(u)\tilde{\gamma}_2 = \gamma_3$ in $P\big(^{\sigma(u)}z,\ ^{\sigma(uv)}z\big)$ under the assumption that $\sigma(u)$ and $\sigma(v)$ commutes with each other. This implies that $Z(uv) = Z(u)Z(v)$ and we are done. ∎

In particular we recover that $Z|_{P_n}$ is a group representation, whose image lies in a copy of $\mathcal{U}(n)$ in $V(B_n)$.

## 4. Braids and regularization of singular integrals.

(4.0) We wish to relate now $Z(uv)$, $Z(u)$ and $Z(v)$ in the general case of not necessary commuting $\sigma(u)$ and $\sigma(v)$. The strategy consists in pushing-off the base points through infinity. Equivalent statements were obtained by Le and Murakami, Bar-Natan, Cartier (see [Le-Mu], [BN2], [Car]). We shall skip over the details which will be considered in the second paper of this series. Our aim is to explain a subtle point around the multiplicativity of iterated integrals and their regularizations.

(4.1) A point is said to be at infinity if it sits on $\mathbb{C}^n \smallsetminus X_n = \Delta_n$. This corresponds to a $n$-string whose points become close to each other depending on the strata of $\Delta_n$ where the limit sits. Incorporating points at infinity would have the effect of replacing the various base points $S_n z$ by an unique base point $z_0$ lying in $\{z_1 = z_2 = \cdots = z_n\} \subset \mathbb{C}^n$. This of course implies the iterated integrals computing the holonomy become singular and regularizations are needed. Such a singular curve has the shape pictured below and usual regularization are those

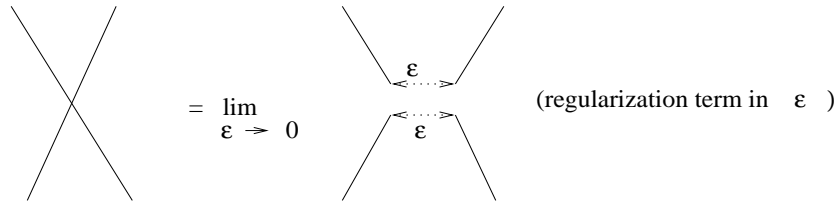

provided by an $\varepsilon$-approximation (see the figure). Now the planar picture of the strings contains always some additional information: the $j^{\text{th}}$ string corresponds to the $j^{\text{th}}$ coordinate $z_j \in \mathbb{C}$. But the role of singularities is just the interchange of two strings. This means that the trajectories in $\mathbb{C}^n$ are like in the picture bellow

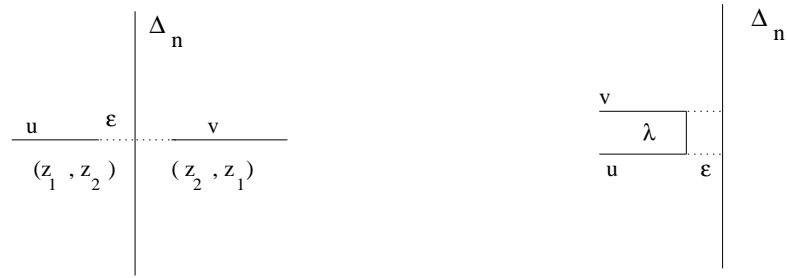

A: real trajectory                          B: what the regularized term computes

As for example in the case $u = v_1$, $v = v_2$ we are computing a limit by inserting some $\varepsilon$-approximations, and the curve corresponding in $\mathbb{C}^n$ has the shape $B$. So that it is not at all clear (but true!) that the two limits $A$ and $B$ when suitably renormalized should be the same; on the right hand side we may apply the multiplicativity of the holonomy to get a closed formula for $Z(uv)$ as $Z(u)Z(\lambda)Z(v)$ (properly renormalized).

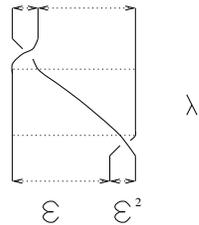

(4.2) To overcome this difficulty it is suitably to work with compactifications of configuration spaces. We recall that, with the notations of §**3**, the general situation is summarized in the picture below and so the endpoints of $^{\sigma(u)}\tilde{\gamma}_2$ and $\gamma_3$

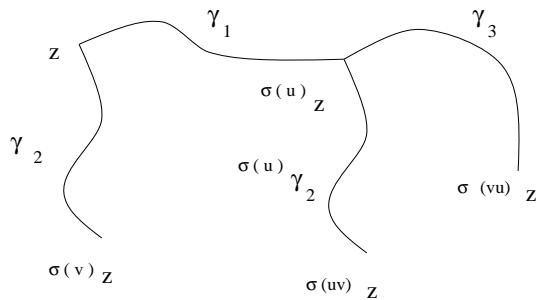

do not coincide. It is for this reason that $z$ is pushed to infinity.

Compactifications of configuration spaces $X_n$ were considered by Fulton and Mac Pherson, Kontsevich and Axelrod and Singer in both algebraic geometric and differential geometric variants. We shall add to $X_n$ a boundary over the $X_k \subset \Delta_n$. We are interested to the strata over the line $\{x_1 = \cdots = x_n\} \subset \Delta_n$; they correspond to binary trees with $n$ labeled leaves. Any real analytic curve $\gamma : [0,1) \to X_n$ which comes from $\gamma : [0,1] \to \mathbb{C}^n$, $\gamma(1) \in \Delta_n$ has a proper lift $\hat{\gamma} : [0,1] \to \widehat{X}_n$ ($\widehat{X}_n$ is the compactification) so that $\hat{\gamma}|_{[0,1)}$ is $\gamma$ when int $\widehat{X}_n$ is identified to $X_n$. Notice that $\hat{\gamma}(1)$ is not uniquely defined by $\gamma(1)$ but also information of how faster $|z_{i+1} - z_i|$ tend to zero is needed. As for example the curves $\gamma_i$ pictured have the endpoints $\hat{\gamma}_i(1)$ for their lifts in the two strata corresponding to the two binary trees of level 2.

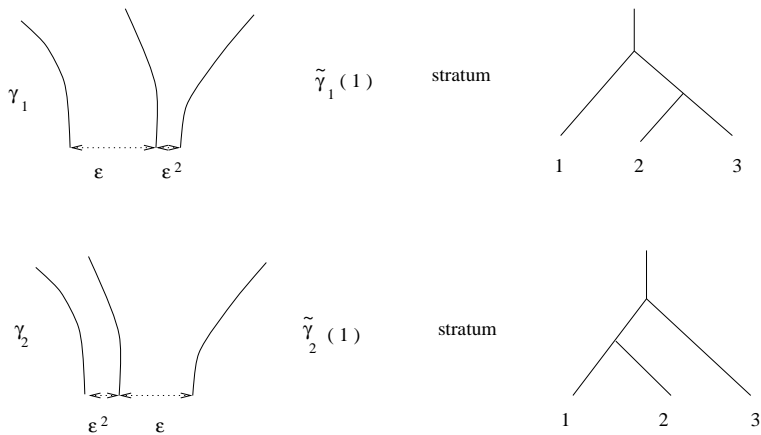

The regularization of the holonomy map is provided by residue theory since the flat connection extends to $\widehat{X}_n$ with regular singularities on the compactification divisors. Now it is easy to see that the singular curves in $X_n$ may have distinct endpoints when lifted to $\widehat{X}_n$, according to the labeling of strings. However we may pass from one stratum to the other using the intermediary curves $\lambda$.

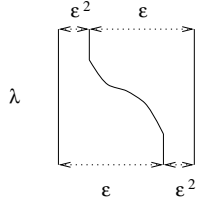

It is for this reason that Le and Murakami's procedure gives the right answer, and the modified (renormalized) integrals provide a representation of their pre-$q$-tangle category.

This type of intermediary curves permit to change the binary trees by fusing moves $F$, and it is easy to see that fusing moves act transitively on the set of binary trees. A more interesting fact is that, once we pass from the tree $T$ to the tree $S$ by a sequence of fusings and compute the holonomy (properly regularized), the result does not depend upon the particular choice of the sequence of fusings we used but only on $S$ and $T$. This is a consequence of the pentagon relation for fusing and is reminiscent to conformal field theory (see [Drin], [BN2]). It follows from the flatness of the connection on $X_n$, and we notice that it is basically the only data which permits to construct link invariants solely (see [BN2]).

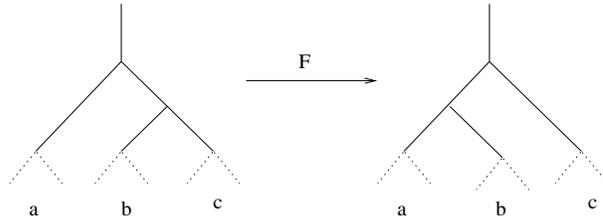

(4.3) Finally the regularization of singular integrals is done in [Le-Mu] and is in some sense canonical: let $\gamma : [0,1) \to X_n$ be a real analytic curve, whose lift $\hat{\gamma}$ has the endpoint $\hat{\gamma}(1)$ in the stratum corresponding to the labeled tree $T$. We represent $\gamma$ as a $n$-string in $C \times \mathbb{R}$.

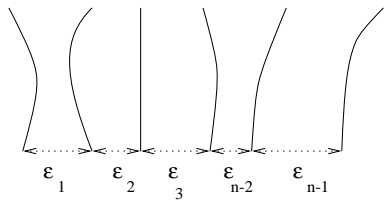

We read from the binary tree $T$ the order in which neighbor points become closer to each other, say $\varepsilon_{\tau(1)}, \ldots, \varepsilon_{\tau(n-1)}$, where $\tau$ is a permutation of $\{1, 2, \ldots, n-1\}$.

We set after [Le-Mu]

$$\varepsilon(T) = \prod_{k=1}^{n-1} \varepsilon_k^{Z(T,k)}$$

$$Z(T,k) = \frac{1}{2\pi i} \sum_{\ell_1 \leq p \leq i_k} \sum_{i_k + 1 \leq q \leq \ell_2} Z_{pq}$$

where
$$\ell_1 = \max\{p;\ p < k \text{ and } \tau^{-1}(p) < k\} + 1,$$
$$\ell_2 = \min\{q;\ q \geq k+1 \text{ and } \tau^{-1}(q) < k\} - 1.$$

All the terms $Z(T, k)$ are commuting with each other. Then the regularization term needed for $Z(\gamma)$ is $\varepsilon(T)$: the limit $\lim_{t \to 1} Z(\gamma|_{[0,t]}) \varepsilon(T)$ is finite, and we denote it by $\widehat{Z}(\gamma)$.

Now if both endpoints of $\gamma$ are at infinity then we have left and right regularization terms as $\lim_{t \to 0} \varepsilon^{-1}(S) Z(\gamma|_{[t, 1-t]}) \varepsilon(T)$ according to the trees $S$, $T$ associated to $\hat{\gamma}(0)$ and $\hat{\gamma}(1)$. These renormalizations are further compatible to the $\lambda$ curves which we insert. The homogeneity of fusing moves permits to compute $\widehat{Z}(\lambda)$ in terms solely of a particular $\lambda$-curve. It is this way Drinfeld's associator $\phi$ appears.

(4.4) We recall from [Drin] that the differential equation
$$G'(x) = \left(\frac{A}{x} + \frac{B}{1-x}\right) G(x), \quad x \in (0, 1)$$
where $G$ is real analytic on $x$, whose coefficients are formal series in two non-commutating variables $A$ and $B$, has unique solutions $G_1$, $G_2$ having prescribed asymptotics
$$G_1(x) \sim x^{A/2\pi i} \quad \text{around } x \sim 0,$$
$$G_2(x) \sim (1-x)^{B/2\pi i} \quad \text{around } x \sim 1.$$
Further $\phi(A, B)$ is defined as the formal series $G_1^{-1}(x) G_2(x)$, and is called Drinfeld's associator. It turns out that the simplest $\lambda$-curve in 3 strings has the regularized holonomy $\phi(Z_{12}, Z_{13})$, just from the definition. In fact the half monodromy, as a function on $\varepsilon$ is a solution of the previous stated equation. Now for a general fusing move $F$ like in the picture, the regularization is $\phi\left(\sum_{\substack{p \in a \\ q \in b}} Z_{pq}, \sum_{\substack{p \in b \\ q \in c}} Z_{p,q}\right)$, where $a, b, c$ are viewed as the sets of labels of leaves issued from the respective vertices.

(4.5) Now in order to find $\widehat{Z}(u)$, $u \in B_n$ we need to fix tree $T_0$ for the initial points. Each $u$ induces a change of the binary tree: in order to be more precise we assume that all the endpoints of the braid drawn in $\mathbb{C} \times \mathbb{R}$ lie on the two lines $\mathbb{R} \times \{0\}$ and $\mathbb{R} \times \{1\}$ so that only the $|z_i - z_{i+1}|$ are taking into account. Then the tree $T_0$ changes into a tree $T_1 = \sigma(u) T_0$, depending only on $\sigma(u) \in S_n$. Let $\phi(u)$ be the product of fusings we need to pass from $T_0$ to $\sigma(u) T_0$. This may be explicitly computed from $\phi$ and $\sigma(u)$, as a product of Drinfeld's associators.

We may state now:

THEOREM 13 (Multiplication law). — *The regularized invariant $\widehat{Z}$ has the following multiplication law*
$$\widehat{Z}(uv) = \widehat{Z}(u) \widehat{Z}(v) \left(\phi([\sigma(u), \sigma(v)], 1)\right)$$
*where $\phi : S_n \to \mathcal{U}(n)$ is the homomorphism defined above.*

Notice we may derive similar relations with $\phi$ inserted between $\widehat{Z}(u)$ and $\widehat{Z}(v)$ or before them.

## 5. Geometric interpretation for $V(n)$.

(5.0) We have a similar result as in Proposition 5, for the algebras $V(n)$. Consider the $S_n$-diagrams constructed like $\mathcal{D}(\underbrace{|\ |\cdots|}_{n})$ but in a more general context:

(1) the vertical arrows can cross each other transversely, this time, and are numbered $1, 2, \ldots, n$;

(2) we have a finite set of horizontal chords whose endpoints are on the vertical arrows;

other;

(4) the isotopy condition: one can move the vertical arrows by preserving the horizontal chord endpoints like in the two moves describes bellow (and coming from the presentation of $S_n$):

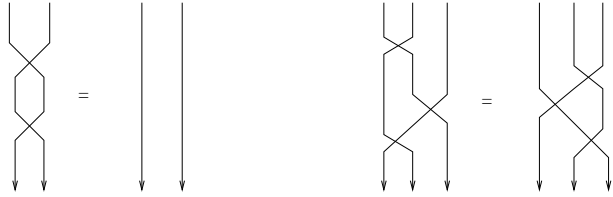

(5) the usual $4T$ relations as in Section 2.

Then the free **k**-module on $S_n$-diagrams quotiened by the equivalence relations (4) and (5) form the algebra of $S_n$-diagrams $TS_n$ extending $J_n$.

PROPOSITION 14. — *The map $V(n) \to TS_n$ given on generators by*

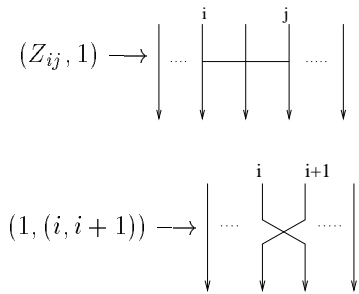

*(where $(i, i+1)$ stands for the transposition interchanging $i$ and $i+1$) is an isomorphism of graded algebras.*

*Proof.* — This map has an obvious inverse, and the $4T$ relations correspond to relations in $\mathcal{U}(n)$, and isotopy moves to the relations in $S_n$. ∎

# APPENDIX

## The product formula

We derive now a direct proof of the product formula of Proposition 9.

Let's look first at $b_1^2 \in P_2$ and consider the representing loop in $X_3$ be that given in the picture below.

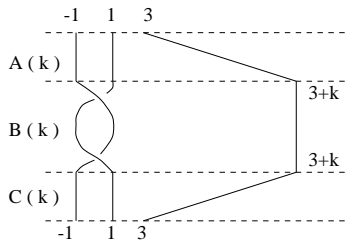

Specifically we have the parametrizations

$$A(k) : z_1 = -1, \ z_2 = 1, \ z_3 = 3 + kt, \ \ t \in [0,1], \ k > 0$$
$$B(k) : z_1 = -\exp(2\pi it), \ z_2 = \exp 2\pi it, \ z_3 = 3 + k, \ \ t \in [0,1]$$
$$C(k) : z_1 = -1, \ z_2 = 1, \ z_3 = 3 + k - kt, \ \ t \in [0,1].$$

We compute the integrals corresponding to $A(k)$, $B(k)$, $C(k)$ separately:

$$Z(B(k)) = \exp(Z_{12}) + Z_{13} \int_0^1 \frac{\exp(2\pi it)}{\exp(2\pi it) + 3 + k} dt + Z_{23} \int_0^1 \frac{\exp(2\pi it)}{\exp(2\pi it) - 3 - k} dt$$
$$+ \sum_{\text{monomials mon},p} \text{mon}(Z_{12}, Z_{13}, Z_{23}) \int_{0 \leq t_1 < \cdots < t_p \leq 1} \prod_{j=1}^\ell \frac{\exp(2\pi it_j)}{\exp(2\pi it_j) \pm (3+k)} dt_1 \cdots dt_p,$$

where the sum is made up on all non-commutative monomials in $Z_{12}$, $Z_{13}$, $Z_{23}$, and $\ell \geq 1$ is the number of $Z_{13}$'s and $Z_{23}$'s in the monomial. We don't need explicit computations of all coefficients. It suffices to see that the coefficient of such a monomial occuring in the sum is bounded as follows:

$$|\text{coefficient}| < \int_{0 \leq t_1 < \cdots < t_p \leq 1} \left| \frac{1}{\exp(2\pi it_j) \pm (3+k)} \right| dt_1 \cdots dt_p < \frac{1}{(k+2)p!}.$$

Next we may compute

$$Z(A(k)) = 1 + Z_{13} \log \frac{4 + kt}{4} + Z_{23} \log \frac{2 + kt}{2}$$
$$+ \sum_{\text{monomials mon},p} \text{mon}(Z_{13}, Z_{23}) \int_{0 \leq t_1 < \cdots < t_p \leq 1} k^p \prod_j \frac{1}{(kt_j + 3) \pm 1} dt_1 \cdots dt_p,$$

where the sum is taken over the monomials in two non-commutative variables of degree at least two. Finally $Z(C(k))$ is the same sum as above but the coefficients of a degree $p$ monomial is multiplied by $(-1)^p$. Again we seek for upper bounds of these coefficients:

$$|\text{coefficient}| < \int_{0 \leq t_1 < \cdots < t_p \leq 1} \frac{k^p}{(kt_p + 2)^p} dt_1 \cdots dt_p \leq Q_p(\log(k+2)),$$

for some polynomial $Q_p$ depending on $p$.

Consequently the product

$$Z(b_1^2) = Z(A(k))Z(B(k))Z(C(k)) = \exp Z_{12} + \text{mon}(Z_{12}, Z_{13}, Z_{23}) \ \text{coeff(mon)},$$

and each coefficient is now

$$|\text{coeff}| \leq \frac{Q_p(\log(k+2))}{k+2}.$$

On the other hand $Z(b_1^2)$ does not depend on the choice of $k$. So letting $k$ goes to infinity gives $Z(b_1^2) = \exp Z_{12}$.

Now the general case goes similarly. For a pure braid $x \in P_n \hookrightarrow P_{n+1}$ we push the $(n+1)^{\text{th}}$ string to infinity as above. Again

$$Z(A(k)) = 1 \sum \text{coeff } (A,k) \ (\text{monomials containing only } Z_{j \ n+1})$$
$$Z(C(k)) = 1 \sum \text{coeff } (C,k) \ (\text{monomials containing only } Z_{j \ n+1})$$
$$Z(X(k)) = Z(x) + \text{coeff } (X,k) \ (\text{new monomials in all } Z_{ij}, \text{ at least}$$
$$\text{one variable being some } Z_{j \ n+1}).$$

and we have estimations

$$|\text{coeff } (A,k)|, \ |\text{coeff } (C,k)| \leq Q_p(\log k),$$

$$|\text{coeff } (X,k)| \leq \frac{1}{R_p(k)},$$

where $Q_p, R_p$ are polynomials depending on the degree $p$ and $x$. We derive that for $k \to \infty$
$$Z(A(k)X(k)C(k)) \longrightarrow Z(x),$$
and the independence of the monodromy on the path yields the result for the inclusion $P_n \hookrightarrow P_{n+1}$. Successive use of this prove the claim.

## References

[Arnold] ARNOLD V.I. — *The cohomology ring of colored braid groups*, Math. Zametki **5** (1969), 227–231.

[BN1] BAR-NATAN D. — *On the Vassiliev knot invariants*, Topology, (to appear).

[BN2] BAR-NATAN D. — *Non-associative tangles*, Proc. Georgia Top. Conf., 1994 (to appear).

[Bir] BIRMAN J. — *New points of view in knot and link theory*, Bull. A.M.S. **28** (1993), 253–287.

[Bir-Lin] BIRMAN J., LIN X.S. — *Knot polynomials and Vassiliev's invariants*, Invent. Math. **111** (1993), 225–270.

[Ber] BERCEANU B. — *Homotopie rationnelle des espaces de petites dimensions*, Thèse Univ. Paris 6, 1994.

[Chen1] CHEN K.T. — *Iterated path integrals*, Bull. A.M.S. **83** (1977), 831–878.

[Chen2] CHEN K.T. — *Extension of $C^\infty$ function algebra by integrals and Malcev completion of $\pi_1$*, Adv. Math. **23** (1977), 181–210.

[Car] CARTIER P. — *Construction combinatoire des invariants de Vassiliev-Kontsevich*, C. R. Acad. Sci. Paris Sér. I Math. **316** (1993), 1205–1210.

[Drin] DRINFELD V. — *On quasi-triangular quasi-Hopf algebras and a group closely connected to $\mathrm{Gal}(\mathbb{Q};\mathbb{Q})$*, Leningrad Math. J. **2** (1991), 829–860.

[DGMS] DELIGNE P., GRIFFITHS P., MORGAN J., SULLIVAN D. — *Real homotopy theory of Kähler manifolds*, Invent. Math. **29** (1975), 245–274.

[FR] FALK M., RANDELL R. — *Pure braid groups and products of free groups*, Contemp. Math. **78** (1988), .

[KV] KASSEL C., TURAEV,V. — *Chord diagram invariants of tangles and graphs*, Preprint Strasbourg **15** (1995), .

[Kohno1] KOHNO T. — *Monodromy representations of braid groups and Yang-Baxter equations*, Ann. Inst. Fourier (Grenoble) **37** (1987), 139–160.

[Kohno2] KOHNO T. — *Série de Poincaré-Koszul associée aux groupes de tresses pures*, Invent. Math. **82** (1985), 57–76.

[Konts] KONTSEVICH M. — *Vassiliev's knot invariants*, Adv. Soviet. Math. **16** (1993), 137–150.

[Le-Mu] LE T.Q.T., MURAKAMI J. — *Representations of the category of tangles by Kontsevich's iterated integrals*, Commun. Math. Phys. **168** (1995), 535–563.

[Piu] PIUNIKHIN S. — *Combinatorial expression for universal Vassiliev link invariant*, Commun. Math. Phys. **168** (1995), 1–22.

[St] STANFORD T. — *The functoriality of Vassiliev-type invariants of links, braids and knotted graphs*, J. Knot Th. and its Ramif. **3** (1994), 247–262.

[Vas] VASSILIEV V.A. — *Cohomology of knot spaces*, Adv. Soviet. Math. **1** (1990), 23–69.

[V] VOGEL P. — *Invariants de Vassiliev des noeuds (d'après D.Bar-Natan, M.Kontsevich et V.A.Vassiliev)*, Sem. Bourbaki **769** (1993), .

– ◇ –

Université de Grenoble I
**Institut Fourier**
Laboratoire de Mathématiques
associé au CNRS (URA 188)
B.P. 74
38402 ST MARTIN D'HÈRES Cedex (France)

(6 octobre 1995)
e-mail: funar@fourier.ujf-grenoble.fr